\documentclass[12pt,twoside]{article}
\usepackage{amsfonts}
\usepackage{amssymb,amsmath}
\begin{document}
\begin{center} \textbf{Finiteness properties of abc-equations c = a+b}\\
\par
\vspace{10pt}
Constantin M. Petridi\\
 cpetridi@math.uoa.gr\\
 cpetridi@hotmail.com
\end{center}
\par
\vspace{5pt}
\begin{center}
\small{
\begin{tabular}{p{11cm}}
\vspace{20pt} \textbf{Abstract.} We classify integer
abc-equations $c=a+b$ (to be defined), according to their radical
R(abc) and prove that the resulting equivalence classes contain
only a finite number of such equations. The proof depends on a
1933 theorem of Kurt Mahler.
\end{tabular}
}
\end{center}
\vspace{60pt}
\par
\textbf{1. abc-equations, their classes and types}
\par
\vspace{20pt} For $i=1,2,\ldots,p_{i},q_{i},r_{i}$ denote primes,
$x_{i},y_{i},z_{i}$ integer variables $\geq 1$.
$R(a)=\prod_{1}^{\omega}p_{i}$ is the radical of
$a=\prod_{1}^{\omega}p_{i}^{x_{i}}$ and $(a,b)$ the g.c.d. of $a$
and $b$. \vspace{20pt}
\par
Integer equations $E:c=a+b$ are defined as abc-equations iff
$1\leq a < b,\;\;(a,b)=1$. The radical $R(E)$ of $E$ is defined
by $R(E)=R(abc)=\prod_{p|abc}\,p=p_{1}p_{2}\cdots p_{\omega}$,
where $p_{1}<p_{2}<\cdots <p_{\omega}$. Clearly $p_{1}=2$.
\vspace{20pt}
\par
We embody all abc-equations with same radical
$R(abc)=p_{1}p_{2}\cdots p_{\omega}$ into one class, denoted
$C_{p_{1}\cdots p_{\omega}}$. This induces a classification of
the totality of all abc-equations into disjoint equivalence
classes $C_{p_{1}\cdots p_{\omega}}$, where $a,b$ run
independently over the $abc$-domain of constraint $$1\leq a <
b,\;\;(a,b)=1.$$
\par
To each abc-equation $E$ belonging to the class $C_{p_{1}\cdots
p_{\omega}}$ we associate its type $T(E)$. To illustrate this
point consider e.g. the equations
$E_{1}:p_{1}^{x_{1}}p_{3}^{x_{2}}=p_{2}^{x_{3}}+p_{4}^{x_{4}}$ and
$E_{2}:p_{1}^{y_{1}}=p_{2}^{y_{2}}p_{4}^{y_{3}}+p_{3}^{y_{4}}$
(assumed satisfiable in $p_{i},x_{i},y_{i},\;\; i = 1,2,3,4$).
They both belong to the class $C_{p_{1}p_{2}p_{3}p_{4}}$, but
there is a formal difference between them. The first is of type
$T(E_{1})=(pp,p,p)$, the second of type $T(E_{2}=(p,pp,p)$. For
the first five cases, $\omega = 2,3,4,5,6$ the possible types are:
\vspace{100pt} \pagebreak
\begin{center}
\begin{tabular}{ccccc}
$\omega=2$&$\omega=3$&$\omega=4$&$\omega=5$&$\omega=6$\\
&&&&\\
$(p,1,p)$&$(p,1,pp)$&$(p,1,ppp)$&$(p,1,pppp)$&$(p,1,ppppp)$\\
         &$(p,p,p)$ &$(p,p,pp)$ &$(p,p,ppp)$ &$(p,p,pppp)$\\
         &$(pp,1,p)$&$(pp,1,pp)$ &$(p,pp,pp)$ &$(p,pp,ppp)$\\
         &          &$(pp,p,p)$ &$(pp,1,ppp)$&$(pp,1,pppp)$\\
         &          &$(ppp,1,p)$&$(pp,p,pp)$ &$(pp,p,ppp)$\\
         &          &           &$(ppp,1,pp)$&$(pp,pp,pp)$\\
         &          &           &$(ppp,p,p)$ &$(ppp,1,ppp)$\\
         &          &           &$(pppp,1,p)$&$(ppp,p,pp)$\\
         &          &           &            &$(pppp,1,pp)$\\
         &          &           &            &$(pppp,p,p)$\\
         &          &           &            &$(ppppp,1,p)$
\end{tabular}
\end{center}
\par
\vspace{15pt}
The type of an abc-equation $c=a+b$ may also be
written as
$$T(E)=(\omega(c),\omega(a),\omega(b)),$$
where $\omega(n)$ is the number of different prime factors of
$n$. Our notation puts more in evidence the combinatorial
peculiarity of this notion.
\par
\vspace{15pt}
 For $\omega \geq 2$ the number of possible types is
$\big[\frac{(\omega+1)^{2}}{4}\big]-1$. We shall not give the
proof (combinatorial) of this here, as for our present purpose it
is only the finiteness of this number that is essential for a
fixed $\omega$.
\par
\vspace{15pt} \textbf{2. Connection with the recurrence
$P(n+\varphi(a))=a+P(n)$}
\par
\vspace{15pt} We refer and use notations and results of our
article [2]. Since in the abc-equation $c=a+b$, $c$ and $b$ are
coprime to $a$, they necessarily do appear in the sequence $P(n)$
with certain indexes. Denoting the index of $b$ by $n$, it follows
that the index of $c$ is $n+\varphi(a)$. All abc-equations,
therefore,
\begin{center}
\begin{tabular}{l}
$c = a+b$\\$1\leq a<b$\\$(a,b)=1,$
\end{tabular}
\end{center}
can also be written as
\begin{center}
\begin{tabular}{c}
$P(n+\varphi(a))=a+P(n)$\\$1\leq a$\\$1\leq n,$
\end{tabular}
\end{center}
where $a$ and $n$ run independently over the indicated domain.
\vspace{15pt}
\par
The radical $R(abc)$ of the abc-equation takes the form
$$R(abc)=R\big\{aP(n)P(n+\varphi(a))\big\}.$$
\par
\vspace{10pt} \textbf{3. Mahler's Theorem}
\par
\vspace{20pt} In his 1933 paper [1] p.724-725, Kurt Mahler proved
following theorem (freely translated from German), as a result of
his investigations on the approximation of algebraic numbers.
\par
\vspace{20pt} Theorem (Mahler). Let $M_{1},\;M_{2},\;M_{3}$ be
finite sets of primes, whose union consists of different primes.
If the prime divisors of the natural number $Z_{1}$ belong to
$M_{1}$, the prime divisors of the natural number $Z_{2}$ to
$M_{2}$ and the prime divisors of the natural number $Z_{3}$ to
$M_{3}$, then the equation
$$Z_{1}+Z_{2}=Z_{3}$$ has only a finite number of solutions.
\vspace{10pt}
\par
In his proof he sets $M_{1}=\{P_{1},P_{2},\ldots,P_{t}\}$,
$P_{i}$ prime, $1\leq i\leq t$,
$M_{2}=\{Q_{1},Q_{2},\ldots,Q_{u}\}$, $Q_{i}$ prime, $1\leq i\leq
u$, $M_{3}=\{R_{1},R_{2},\ldots,R_{v}\}$, $R_{i}$ prime, $1\leq
i\leq v$, so that the equation $Z_{1}+Z_{2}=Z_{3}$ becomes
$$P_{1}^{p_{1}}P_{2}^{p_{2}}\cdots P_{t}^{p_{t}}+Q_{1}^{q_{1}}Q_{2}^{q_{2}}\cdots
Q_{s}^{p_{s}}= R_{1}^{r_{1}}R_{2}^{r_{2}}\cdots R_{v}^{r_{v}},$$
where the exponents
$p_{1},p_{2},\ldots,p_{t},q_{1},q_{2},\ldots,q_{u},r_{1},r_{2},\ldots,r_{v}$
are non-negative integers.
\vspace{15pt}
\par
It is this form of equation we shall use in our application of
Mahler's theorem.
\vspace{40pt}
\par
\textbf{4. The finiteness of the classes $C_{p_{1}\cdots
p_{\omega}}$}
\par
\vspace{20pt} \textbf{Theorem.} For a given fixed radical
$p_{1}\cdots p_{\omega}$ the class $C_{p_{1}\cdots p_{\omega}}$
contains only a finite number of abc-equations.
\par
\vspace{5pt}\textbf{Proof.} Suppose the opposite is true. Then
there is an infinite number of abc-equations $c=a+b$ with radical
$R(abc)=p_{1}\cdots p_{\omega}$.
\vspace{10pt}
\par
Since, according to $\S2$, for a given fixed $\omega$ there are
only finitely many types, there must be at least one type
$T(\,\underbrace{p\cdots p}_{\kappa}\,,\underbrace{p\cdots
p}_{\lambda}\,, \underbrace{p\cdots p}_{\mu}\,)$,
$\kappa+\lambda+\mu=\omega$ for which there is an infinite number
of abc-equations with $R(abc)=p_{1}\cdots p_{\omega}$. Let these
equations be \pagebreak
$$q_{1}^{x_{i1}}\cdots q_{\kappa}^{x_{i\kappa}}\;=\;
r_{1}^{y_{i1}}\cdots r_{\lambda}^{y_{i\lambda}}\;+\;
s_{1}^{z_{i1}}\cdots s_{\mu}^{z_{i\mu}}$$
$$x_{ij}\geq 1,\;\;1\leq j\leq \kappa$$
$$y_{ij}\geq 1,\;\;1\leq j\leq \lambda$$
$$z_{ij}\geq 1,\;\;1\leq j\leq \mu$$
$$i = 1,2,\ldots,\infty$$
\par where $\{q_{1},\ldots
q_{\kappa}\},\;\{r_{1},\ldots r_{\lambda}\},\;\;\{s_{1},\ldots
s_{\mu}\}$ are disjoint sets of primes satisfying
$$\{q_{1},\ldots q_{\kappa}\}\;\cup\; \{r_{1},\ldots
r_{\lambda}\}\;\cup\;\{s_{1},\ldots s_{\mu}\}\;=\;\{p_{1},\ldots
,p_{\omega}\}.$$
\par Above equations with their respective constraints are exactly
those considered by Mahler, excluding the equations where at
least one of the exponents
$p_{1},p_{2},\ldots,p_{t},q_{1},q_{2},\ldots,q_{u},r_{1},r_{2},\ldots,r_{v}$
is zero. Applying his theorem it follows that the supposition
that there are infinite such equations leads to a contradiction.
\vspace{10pt}
\par
Hence the class $C_{p_{1}\cdots p_{\omega}}$
contains only a finite number of abc-equations. $\Box$
\par
\vspace{15pt}\textbf{Note.} Above result is an immediate
consequence of the abc-conjecture, assumed to be true.
\par
\vspace{20pt} \textbf{Acknowledgment.} I am indebted to Peter
Krikelis, University of Athens, Department of Mathematics, for his
 assistance.
\par
\vspace{20pt} \textbf{References}\vspace{10pt}
\par
[1] Kurt Mahler, Zur Approximation algebreischer Zahlen. I.
(\"{U}ber den gr\"{o}ssten Printeiler bin\"{a}rer Formen), Math.
Ann. 107(1933),p.691-730.
\par
\vspace{10pt}[2] Constantin M. Petridi, The integer recurrence
$P(n) = a+P(n-\varphi(a))$ I, arXiv:1208.5348v2 [math.NT] 30 Aug
2012

\end{document}